\documentclass[11pt]{amsart}

\usepackage{url}
\usepackage{hyperref}
\usepackage{amssymb}
\usepackage{amsmath}
\usepackage{mathrsfs}
\input xy
\xyoption{all}
\usepackage{mathdots}
\usepackage{mathtools}
\usepackage{youngtab}
\usepackage{enumerate}
\usepackage{graphicx}  
\usepackage{tikz}
\usepackage{tkz-graph}

\newtheorem{thm}{Theorem}

\newtheorem{cor}[thm]{Corollary}
\newtheorem{prop}[thm]{Proposition}

\newtheorem{Definition}[thm]{Definition}
\newenvironment{definition}
  {\begin{Definition}\rm}{\end{Definition}}

\newtheorem{Example}[thm]{Example}
\newenvironment{example}
  {\begin{Example}\rm}{\end{Example}}
  
\newtheorem{Remark}[thm]{Remark}
\newenvironment{remark}
  {\begin{Remark}\rm}{\end{Remark}}

\title{RSK via local transformations}
\author{Sam Hopkins}
\date{July, 2014}

\begin{document}

\maketitle

%IMPROVE INTRODUCTION, OBV, and relocate the remark at the end

\section{Introduction}

The Robinson--Schensted--Knuth correspondence is a map from~\mbox{$n \times n$ $\mathbb{N}$}-matrices (i.e., matrices with nonnegative integer entries) to pairs of semi-standard Young tableaux of the same shape. For basic definitions of and background on tableaux and the RSK algorithm, see~\cite[Ch. 7]{stanley2}. Abstractly, we may regard RSK as a bijection
\begin{align*}
\{ \textrm{$n \times n$ $\mathbb{N}$-matrices} \} &\xrightarrow{RSK} \left \{ \parbox{3in}{ \begin{center} $(P,Q)\colon P,Q \textrm{ are SSYT with } \mathrm{sh}(P) = \mathrm{sh}(Q)$ \\ \textrm{and whose entries are in $[n] := \{1,\ldots,n
\}$} \end{center} } \right \} \\
A &\xmapsto{RSK} (P,Q)
\end{align*}
with the property that $\mathrm{type}(P) = \mathrm{col}(A)$ and $\mathrm{type}(Q) = \mathrm{row}(A)$. Here for a SSYT $T$ with entries in $[n]$, the \emph{type} of $T$ is defined to be
\[ \mathrm{type}(T) := (\mathrm{type}(T)_1,\ldots,\mathrm{type}(T)_n),\]
where 
\[\mathrm{type}(T)_i := \textrm{the number of entries of $T$ equal to $i$},\] 
and for a $n \times n$ $\mathbb{N}$-matrix~$B = (b_{ij})$ the \emph{column} and \emph{row sums} of $B$ are
\begin{align*} 
\mathrm{col}(B) &:= (\mathrm{col}(B)_1, \ldots, \mathrm{col}(B)_n) := (\sum_{i=1}^{n} b_{i1}, \ldots, \sum_{i=1}^{n} b_{in}) \\
\mathrm{row}(B) &:= (\mathrm{row}(B)_1, \ldots, \mathrm{row}(B)_n) := (\sum_{i=1}^{n} b_{1i}, \dots, \sum_{i=1}^{n} b_{ni}).
\end{align*}
RSK is normally defined in terms of a somewhat complicated procedure called row-insertion. This definition unfortunately masks the remarkable symmetry property that if $A \xmapsto{RSK} (P,Q)$ then $A^t \xmapsto{RSK} (Q,P)$. Here we present an alternate description of RSK via local transformations that makes the symmetry obvious. Our approach is thus similar in spirit to the growth diagrams of Fomin~\cite{fomin}. In order to motivate this description, we first examine the Gelfand-Tsetlin patterns associated to the Young tableaux that result from RSK, and show how RSK can be thought of as a map between matrices.

\begin{definition}
A \emph{Gelfand-Tsetlin pattern} is a triangular array of nonnegative integers 
\[\begin{array}{c}
g_{1,1} \qquad \qquad g_{1,2} \qquad \qquad g_{1,3} \qquad \qquad  \cdots \qquad  \qquad g_{1,n} \\
g_{2,2} \qquad \qquad  g_{2,3} \qquad \qquad \cdots \qquad \qquad g_{2,n} \\
g_{3,3} \qquad \qquad  \cdots \qquad \qquad g_{3,n} \\
\ddots \qquad \qquad \qquad \iddots \\
g_{n,n}
\end{array}\]
such that $g_{i,j} \geq g_{i+1,j+1} \geq g_{i,j+1}$ for all $1\leq i \leq j \leq n-1$. Gelfand-Tsetlin patterns with $n$ rows are in bijection with semi-standard Young tableaux with entries among $[n]$ in a natural way, as explained below. (There seems to be no consensus in the literature about how Gelfand-Tsetlin patterns are indexed so we are merely choosing one convention here.)
\end{definition}

Suppose $A \xmapsto{RSK} (P,Q)$. We associate an $n$-row Gelfand-Tsetlin pattern to $P$, denoted $GT(P)$, as follows. Define the series $P^{1},\ldots,P^{n}$ of SSYT by
\begin{itemize}
\item $P^{1} := P$;
\item $P^{i+1}$ is obtained from $P^{i}$ by removing all entries of $n+1-i$.
\end{itemize}
Then $GT(P)$ is just the Gelfand-Tsetlin pattern whose $i$-th row is $\mathrm{sh}(P^{i})$. We define $GT(Q)$ similarly. In other words, for a SSYT $T$ with entries in~$[n]$, we define $GT(T)$ to be the Gelfand-Tsetlin pattern with $n$ rows defined by
\[ g_{i,j} := \textrm{number of entries $\leq n+1-i$ in row $j-i+1$ of $T$}.\]
It is easy to check that the necessary inequalities on the entries hold. Note also that because $\mathrm{sh}(P) = \mathrm{sh}(Q)$, the first rows of~$GT(P)$ and $GT(Q)$ are the same. So it is possible to glue $GT(P)$ and~$GT(Q)$ along their first row to form a single $n\times n$ $\mathbb{N}$-matrix that is weakly increasing in rows and columns. Formally, suppose $GT(P) = (g_{ij})$ and $GT(Q) = (h_{ij})$; then the matrix $\widehat{A} = (\widehat{a}_{ij})$ is given by
\[ \widehat{a}_{ij} := \left\{
	\begin{array}{ll}
		g_{i-j+1,n+1-j}  & \textrm{if $i \geq j$;} \\
		h_{j-i+1,n+1-i} & \textrm{otherwise}.
	\end{array} \right.\] 
This series of combinatorial objects is best understood through an example.

\begin{example} \label{ex:rsk}
Suppose
\[ A = \begin{pmatrix} 1 & 0 & 2 \\ 0 & 2 & 0 \\ 1 & 1 & 0 \end{pmatrix}.\]
Then
\[ A \xmapsto{RSK} \Bigg( P = \raisebox{-.25in}{\young(1122,23,3)}, Q =  \raisebox{-.25in}{\young(1113,22,3)} \; \Bigg),\]
and
\begin{align*}
GT(P) &= {\begin{array}{c}
4 \qquad  2 \qquad 1 \\
4 \qquad 1\\
2
\end{array}} \\ \\
GT(Q) &= {\begin{array}{c}
4 \qquad 2 \qquad 1 \\
3 \qquad 2\\
3
\end{array}}
\end{align*}
so
\[ \widehat{A} = \begin{pmatrix} 1 & 2 & 3 \\ 1 & 2 & 3 \\ 2 & 4 & 4 \end{pmatrix}.\]
Indeed, $\widehat{A}$ is weakly increasing along rows and columns. $\square$
\end{example}

We have built up a series of maps between combinatorial objects
\[ A \mapsto (P,Q) \mapsto (GT(P),GT(Q)) \mapsto \widehat{A}.\]
By ignoring the intermediary steps, we can really think of RSK as a map
\begin{align*}
\{ \textrm{$n \times n$ $\mathbb{N}$-matrices} \} &\xrightarrow{RSK} \left \{ \parbox{2.25in}{\begin{center} \textrm{$n \times n$ $\mathbb{N}$-matrices that are weakly}  \\ \textrm{increasing in rows and columns} \end{center}} \right \} \\
A &\xmapsto{RSK} \widehat{A}.
\end{align*}
Because each intermediary step was bijective, it is clear that this map is a bijection. There are a few more properties of this map worth noting. First of all, the diagonal sums of $\widehat{A}$ record information about the row and column sums of $A$.

\begin{definition}
For an $n \times n$ $\mathbb{N}$-matrix $B = (b_{ij})$ and some $i,j \in [n]$, define the \emph{diagonal sum at entry $(i,j)$}, denoted $\mathrm{diag}_B(i,j)$, to be
\[ \mathrm{diag}_B(i,j) := \sum_{k = 0}^{\mathrm{min}(i,j)-1} b_{i-k,j-k}.\]
\end{definition}

Using the above notation, when $A \xmapsto{RSK} \widehat{A}$, we have for any $i \in [n]$
\begin{align*}
\mathrm{diag}_{\widehat{A}}(n,i) &= \sum_{k=1}^{i} \mathrm{col}(A)_k; \\
\mathrm{diag}_{\widehat{A}}(i,n) &= \sum_{k=1}^{i} \mathrm{row}(A)_k.
\end{align*}
To see why this is the case, observe that $\mathrm{diag}_{\widehat{A}}(n,i)$ is just the sum of the entries in row $(n+1)-i$ of $GT(P)$. But the sum of this row of $GT(P)$ is the number of $1$'s, $2$'s, ..., and $i$'s in $P$. But this is just $\sum_{k=1}^{i} \mathrm{type}(P)_k$, and because $\mathrm{type}(P) = \mathrm{col}(A)$, the diagonal sum is as claimed. Similarly for~$Q$ and row sums. The final key property of our map is that if $A \mapsto \widehat{A}$ then~$A^t \mapsto \widehat{A}^t$, which follows from the above-mentioned symmetry property of the map to Young tableaux.

To summarize, RSK is secretly a map $A \xmapsto{RSK} \widehat{A}$ from $n \times n$ $\mathbb{N}$-matrices to $n \times n$ $\mathbb{N}$-matrices that are weakly increasing in rows and columns such that
\begin{enumerate}[(a)]
\item the map is bijective;
\item $\mathrm{diag}_{\widehat{A}}(n,i) = \sum_{k=0}^{i} \mathrm{col}(A)_k$ and $\mathrm{diag}_{\widehat{A}}(i,n) = \sum_{k=0}^{i} \mathrm{row}(A)_k$ for $i \in [n]$;
\item $A^t \xmapsto{RSK} \widehat{A}^t$.
\end{enumerate}
Our aim is to give an explicit construction of such a map where these essential properties are immediate from the construction.

\begin{remark}
These notes are not intended for publication. They record some ideas that were explained by Alex Postnikov during the MIT combinatorics preseminar. The ideas presented here are closely related to a certain birational lifting of RSK originally studied by Kirillov~\cite{kirillov} which takes matrices with entries in $\mathbb{R}_{> 0}$ to other matrices of this form (or more generally, to certain three-dimensional arrays with entries in $\mathbb{R}_{> 0}$). Birational RSK has recently seen significant interest; see, e.g.,~\cite{noumi},~\cite{danilov1},~\cite{danilov3},~\cite{corwin},~\cite{oconnell1},~\cite{oconnell2}. The work of Danilov and Koshevoy~\cite{danilov1}~\cite{danilov2} in particular is closely related to the presentation of RSK given here because they explicitly observe the connection between RSK and the octahedron recurrence. Although everything written here naturally generalizes to the birational setting we choose to work in the combinatorial max-plus semiring setting because we are just trying to give an alternative description of classical RSK. In particular, our aim is to show how many classical combinatorial results (for instance, the ``hook-length formula'') follow easily from this treatment of RSK, and at a greater level of generality as well. Thanks to Darij Grinberg for proofreading and helpful comments.
\end{remark}

\begin{remark}
In May of 2017, Igor Pak pointed out his paper~\cite{pak} to me, in which he develops essentially the same piecewise-linear ``toggle'' description of RSK that we will explain in the next section. The toggles can be seen specifically in~\S4 of~\cite{pak}. That paper also cites an earlier paper of Berenstein and Kirillov~\cite{berenstein}, who independently discovered a similar toggle map (and in fact, in~\cite{berenstein2} Berenstein--Kirillov studied the toggle description of RSK).
\end{remark}

\section{The RSK algorithm via toggles}

In fact, our map will work in a slightly more general setting than matrices.

\begin{definition}
An \emph{$\mathbb{N}$-tableau of shape $\lambda$} is an assignment of nonnegative integers to the boxes of some partition $\lambda$. A \emph{(weak) reverse plane partition of shape $\lambda$} is an $\mathbb{N}$-tableau of shape $\lambda$ that is weakly increasing in rows and columns. We will always use matrix coordinates for the boxes of partitions (which, as is already evident, we write in English notation). This means that the upper-left box, which belongs to every partition except the empty partition, is $(1,1)$.
\end{definition}

The map we are about to describe takes $\mathbb{N}$-tableaux of shape $\lambda$ to (weak) reverse plane partitions of shape~$\lambda$. We will write \mbox{$T \xmapsto{\mathcal{RSK}} \widehat{T}$} with calligraphic script to differentiate this procedure from classical RSK. These reverse plane partitions are termed \emph{weak} because we allow $0$ as the value of a box, but we will drop this adjective from now on. Items (a) and~(c) in the above list of desiderata are clear enough in this context, but we need to explain what exactly (b) should mean. We need a little terminology for this purpose

\begin{definition} For any partition $\lambda$, a \emph{border box} of $\lambda$ is a box~$(i,j)$ such that~$(i+1,j+1)$ is not a box. A \emph{corner box} of $\lambda$ is a box $(i,j)$ such that both~$(i+1,j)$ and $(i,j+1)$ are not boxes. When we speak of border or corner boxes of matrices, we mean when they are considered as square partitions. Given an $\mathbb{N}$-tableau $T = (t_{kl})$ of shape $\lambda$, we define the \emph{row} and \emph{column} sums exactly as above for matrices. For some box $(i,j)$ of $\lambda$, we define the \emph{diagonal sum} at $(i,j)$ as above for matrices. Finally, we define the \emph{rectangle sum} at $(i,j)$ to be
\[\mathrm{rect}_T(i,j) := \sum_{k = 1}^{i} \sum_{l = 1}^{j} t_{kl}.\]
Note that if $(i,j)$ is a border box (which is the case that will concern us most often) then
\[\mathrm{rect}_T(i,j) = \sum_{k=1}^{j} \mathrm{col}(T)_k - \sum_{k=1}^{n-i} \mathrm{row}(T)_{(n+1)-k}.\]
\end{definition}

For any matrix $A$ the sum of all the row sums equals the sum of all the column sums. So we can think of (b) as saying that the diagonal sum of $\widehat{A}$ at some border box~$(i,j)$ of the matrix can be obtained by traveling through border boxes from $(n,0)$ to $(i,j)$, where each time you step rightwards you add the corresponding column sum, and each time you step upwards you subtract the corresponding row sum. (We start at $(n,0)$, which is not actually a box of $A$, because the first step rightward to $(n,1)$ corresponds to adding the first column sum.) That is, statement (b) is equivalent to
\[ \mathrm{diag}_{\widehat{A}}(i,j) = \sum_{k=1}^{j} \mathrm{col}(A)_k - \sum_{k=1}^{n-i} \mathrm{row}(A)_{(n+1)-k} = \mathrm{rect}_{A}(i,j)\]
for every $i,j \in [n]$ with at least one of $i$ or $j$ equal to $n$. A good sanity check is to verify that $\widehat{A}$ has the claimed diagonal sums in Example~\ref{ex:rsk}. The translation of (b) for $\mathbb{N}$-tableaux $T$ is thus: for any border box $(i,j)$ of~$\lambda := \mathrm{sh}(T)$ we need $\mathrm{diag}_{\widehat{T}}(i,j) = \mathrm{rect}_T(i,j)$.

We are prepared to give our account of $\mathcal{RSK}$. The idea is to work box-by-box. Given an $\mathbb{N}$-tableau $T'$ we will construct $\widehat{T'}$ with $T' \xmapsto{\mathcal{RSK}} \widehat{T'}$. If~$T' = \emptyset$ then $\widehat{T'} := \emptyset$. So say $T'$ is obtained from~$T$ by adding a single corner box~$(i,j)$ whose entry is $x$. And say~$T \xmapsto{\mathcal{RSK}} \widehat{T}$. Then from $\widehat{T}$ we obtain~$\widehat{T'}$ as follows. Let $\beta_k$ be the entry of $\widehat{T}$ at position $(i-k,j-k)$, let $\gamma_k$ be the entry at position $(i-k+1,j-k)$, and let $\alpha_k$ be the entry at~$(i-k,j-k+1)$ for all~$k \geq 1$. If any of these boxes do not exist in $\widehat{T}$ we regard these values as~$0$; so eventually they are all $0$ for the nonpositive positions. These entries occupy three consecutive diagonals in $\widehat{T}$, as Figure~\ref{fig:diags} illustrates. Our procedure only modifies the values of the boxes that lie along the diagonal where the new box is to be inserted (which were previously the~$\beta_k$). Observe that because~$\widehat{T}$ is weakly increasing in rows and columns, we have for all $k$ the inequality 
\[\mathrm{max}(\alpha_{k+1},\gamma_{k+1}) \leq \beta_k \leq \mathrm{min}(\alpha_k, \gamma_k).\] 
The trick is to ``toggle'' $\beta_k$ between these two extrema. That is, we set in~$\widehat{T'}$ the value of the box previously occupied by $\beta_k$ to be
\[ \beta'_k := \mathrm{max}(\alpha_{k+1},\gamma_{k+1}) + \mathrm{min}(\alpha_k, \gamma_k) - \beta_k,\]
for all $1 \leq k < \mathrm{min}(i,j)$. Finally note that the value of box $(i,j)$ must be greater than or equal to $\mathrm{max}(\alpha_1, \gamma_1)$; so we set it equal to $\mathrm{max}(\alpha_1, \gamma_1) + x$.

\begin{figure}
\newcommand{\bone}{\beta_1}
\newcommand{\btwo}{\beta_2}
\newcommand{\gone}{\gamma_1}
\newcommand{\gtwo}{\gamma_2}
\newcommand{\aone}{\alpha_1}
\newcommand{\atwo}{\alpha_2}
\newcommand{\dotsy}{\scalebox{0.6}{$\ddots$}}
\newcommand{\blank}{\;}
\newcommand{\qwest}{*}
\begin{center}
\scalebox{1.5}{\young(\dotsy \dotsy\blank\blank\blank,\dotsy\btwo\atwo\blank\blank,\blank\gtwo\bone\aone,\blank\blank\gone\qwest,\blank\blank)} \\

\end{center}
\caption{The $\mathbb{N}$-tableau $\widehat{T}$, where position $(i,j)$ has been marked by a star.}
\label{fig:diags}
\end{figure}

\begin{example}
Before we prove the correctness of this algorithm, we give a few simple examples of how it behaves. First of all, for the simplest case of a line $\mathcal{RSK}$ merely records partial sums:
\[ \begin{tabular}{| c | c | c | c |}
\hline
$a_1$ & $a_2$ & $a_3$ & $\ldots$ \\
\hline
\end{tabular} \xmapsto{\mathcal{RSK}} \begin{tabular}{| c | c | c | c |} \hline
$a_1$ & $a_1 + a_2$ & $a_1 + a_2 + a_3$ & $\ldots$ \\ \hline
\end{tabular} \]
Similarly for a hook:
\[ \begin{tabular}{| c | c | c | c |} \hline
$a$ & $b_1$ & $b_2$ & $\ldots$ \\ \hline
$c_1$ \\ \cline{1-1}
$c_2$ \\ \cline{1-1}
$\vdots$ \\ \cline{1-1}
\end{tabular} \xmapsto{\mathcal{RSK}} \begin{tabular}{| c | c | c | c |} \hline
$a$ & $a + b_1$ & $a + b_1 + b_2$ & $\ldots$ \\ \hline
$a + c_1$ \\ \cline{1-1}
$a + c_1 + c_2$ \\ \cline{1-1}
$\vdots$ \\ \cline{1-1}
\end{tabular} \]
Let us consider the first non-hook shape, a $2 \times 2$ square. We already know that
\[ \begin{tabular}{| c | c |} \hline
$a$ & $b$  \\ \hline
$c$ \\ \cline{1-1}
\end{tabular} \xmapsto{\mathcal{RSK}} \begin{tabular}{| c | c |} \hline
$a$ & $a+b$  \\ \hline
$a+c$ \\ \cline{1-1} \end{tabular} \]
So say we add a box in $(2,2)$ with entry $d$. We then must ``toggle'' the value of $(1,1)$; its previous value was $a$ and it is subject to the inequality
\[ 0 \leq a \leq \mathrm{min}(a+b,a+c),\]
so we set it equal to $\mathrm{min}(a+b,a+c) - a = \mathrm{min}(b,c)$. Finally, we set the value of $(2,2)$ to be $\mathrm{max}(a+b,a+c) + d = a + \mathrm{max}(b,c) + d$. Thus
\[ \begin{tabular}{| c | c |} \hline
$a$ & $b$  \\ \hline
$c$ & $d$ \\ \hline
\end{tabular} \xmapsto{\mathcal{RSK}} \begin{tabular}{| c | c |} \hline
$\mathrm{min}(b,c)$ & $a+b$  \\ \hline
$a+c$ & $a + \mathrm{max}(b,c) + d$ \\\hline \end{tabular} \]
It is easily verified that the diagonal sums are correct. $\square$
\end{example}

\begin{prop} \label{prop:welldefined}
The $\mathcal{RSK}$ map defined above is well-defined; that is, the order in which we decompose an $\mathbb{N}$-tableaux~$T$ does not matter.
\end{prop}

\noindent {\bf Proof}:  Suppose $\mathrm{sh}(T) = \lambda$ and let $P_{\lambda}$ be the poset on the boxes of $\lambda$ whereby $(i,j) \leq (k,l)$ iff $i \leq k$ and $j \leq l$. So the upper-left box is minimal in $P_{\lambda}$. The $\mathcal{RSK}$ algorithm described above proceeds by adding boxes one at a time, starting from the empty tableau, according to any order that is a linear extension of $P_{\lambda}$. We need to show that all choices of linear extensions of $P_{\lambda}$ give the same result. As any linear extension of a finite poset can be obtained from any other via a series of adjacent transpositions of incomparable elements, it suffices to show that if $T''$ is obtained from~$T'$ by adding entries in positions $(i,j)$ and~$(i',j')$, both of which are corner boxes of $T''$, then the order in which we add these two entries does not affect how $\mathcal{RSK}$ acts on~$T''$. But note that because $(i,j)$ and~$(i',j')$ are both corner boxes, the diagonals on which $(i,j)$ and  $(i',j')$ lie cannot be adjacent. Therefore, the values of the diagonals adjacent to the diagonal on which~$(i',j')$ lies are not modified when we add some entry at~$(i,j)$ and carry out the local transformations, and vice-versa. But these local transformations \emph{only depend on adjacent diagonals}. So indeed the order in which we add the two entries does not matter. $\square$

\begin{prop}
$\mathcal{RSK}$ is a bijection between $\mathbb{N}$-tableaux of shape $\lambda$ and reverse plane partitions of shape $\lambda$.
\end{prop}

\noindent {\bf Proof}: We can explicitly give the inverse map $\widehat{T'} \xmapsto{\mathcal{RSK}^{-1}} T'$. The crucial fact is that the ``toggle'' steps are reversible. If $\widehat{T'}$ is empty, then so is~$T'$. So consider any corner box $(i,j)$ of $\widehat{T'}$. Let~$\widehat{x}$ be the value of~$\widehat{T'}$ in $(i,j)$; then define $x$ to be the $\widehat{x}$ minus the maximum of the values of $\widehat{T'}$ in the two boxes~$(i-1,j)$ and $(i,j-1)$ (where these are $0$ if there are no such boxes). Let~$\widehat{T}$ be obtained from $\widehat{T'}$ by removing $(i,j)$ and re-toggling along its diagonal. Suppose $\widehat{T} \xmapsto{\mathcal{RSK}^{-1}} T$. Then define $T'$ to be the tableau obtained from $T$ by adding the entry $x$ to position $(i,j)$. This procedure is easily seen to locally reverse the $\mathcal{RSK}$ procedure, so it indeed gives the inverse. $\square$

%start updating here

\begin{prop} \label{prop:diag_rect_sums}
If $T \xmapsto{\mathcal{RSK}} \widehat{T}$, then for any border box $(a,b)$ of $\mathrm{sh}(T)$ we have $\mathrm{diag}_{\widehat{T}}(a,b) = \mathrm{rect}_{T}(a,b)$.
\label{prop:diagsum}
\end{prop}

\noindent {\bf Proof}: Again we consider the local transformations. Suppose $T' \xmapsto{\mathcal{RSK}} \widehat{T'}$, with $T'$ obtained from $T$ by adding an entry of $x$ in a corner box $(i,j)$ and~$T \xmapsto{\mathcal{RSK}} \widehat{T}$. We may assume inductively that the result holds for $T$ and $\widehat{T}$. First observe that the right-hand side of this equality only changes for $(i,j)$ because for any other border box $(i',j')$ of $\mathrm{sh}(T')$, either $j' \geq j$ and $i' < i$ and the contribution of $x$ cancels between $\mathrm{col}(T)_j$ and $-\mathrm{row}(T)_i$, or $j' < j$ and~$i' \geq i$ and there is no contribution of $x$ at all. But observe further that the left-hand side of the equality also only changes for $(i,j)$, because this is the only diagonal whose values we have modified. Thus we need only check that the equality holds at $(i,j)$. Let $\alpha_k,\beta_k,\gamma_k$ be defined as in Figure~\ref{fig:diags}. Then
\begin{align*}
\mathrm{diag}_{\widehat{T'}}(i,j) &= x + \sum_{k\geq 1}(\mathrm{max}(\alpha_k,\gamma_k) + \mathrm{min}(\alpha_k,\gamma_k) - \beta_k) \\
&= x + \sum_{k \geq 1}( \alpha_k + \gamma_k - \beta_k )\\
&= x + \mathrm{diag}_{\widehat{T}}(i,j-1) + \mathrm{diag}_{\widehat{T}}(i-1,j) - \mathrm{diag}_{\widehat{T}}(i-1,j-1).
\end{align*}
But observe that 
\[x + \mathrm{diag}_{\widehat{T}}(i-1,j) - \mathrm{diag}_{\widehat{T}}(i-1,j-1) = x + \mathrm{col}(T)_j = \mathrm{col}(T')_j,\]
by our induction hypothesis, so
\[\mathrm{diag}_{\widehat{T'}}(i,j) = \mathrm{diag}_{\widehat{T}}(i,j-1) + \mathrm{col}(T')_j, \]
and by induction indeed this diagonal sum is as claimed. $\square$

\begin{prop}
If $T \xmapsto{\mathcal{RSK}} \widehat{T}$, then $T^t \xmapsto{\mathcal{RSK}} \widehat{T}^t$.
\end{prop}

\noindent {\bf Proof}: The local transformations evidently commute with transposition, so the whole map does as well. $\square$

\section{An application: hook-length formulas} \label{sec:hook}

We have succeeded in describing via simple local transformations a map with the desired properties (a)--(c). In fact, one further nice property comes along for free. Recall that for a partition $\lambda$, the \emph{hook} of $\lambda$ at $(i,j)$, which we denote $H_{\lambda}(i,j)$, is the set of all boxes weakly below or to the right of $(i,j)$, and the \emph{hook-length} is given by $h_{\lambda}(i,j) := |H_{\lambda}(i,j)|$. (Note that if $(i,j)$ is not a box of $\lambda$, then $H_{\lambda}(i,j) = \emptyset$ and $h_{\lambda}(i,j) = 0$.) And recall that the \emph{weight} of a reverse plane partition $T$, denoted $|T|$, is the sum of all of its entries.

\begin{prop}
Suppose $T \xmapsto{\mathcal{RSK}} \widehat{T}$ with $T = (t_{ij})$ and $\mathrm{sh}(T) = \lambda$. Then
\[ |\widehat{T}| = \sum_{(i,j)\in\lambda} t_{ij}h_{\mathrm{sh}(T)}(i,j). \]
\end{prop}
\noindent {\bf Proof}: One way to write~$|\widehat{T}|$ is as the sum over all border boxes $(a,b)$ of $\mathrm{diag}_{\widehat{T}}(a,b)$. And the number of border boxes $(a,b)$ for which a given $t_{ij}$ contributes to $\mathrm{rect}_{T}(a,b)$ is exactly $h_{\mathrm{sh}(T)}(i,j)$ (these are precisely the border boxes which lie ``between'' the extreme boxes of $H_{\lambda}(i,j)$). So this proposition follows from Proposition~\ref{prop:diag_rect_sums}. $\square$

One easily obtains as a corollary of this last proposition the generating function for reverse plane partitions in terms of hook-lengths (see~\cite{hillman} or~\cite[\S7.22]{stanley2}):

\begin{cor}
Let $\lambda$ be a partition. Then
\[ \sum_{T} q^{|T|} = \prod_{(i,j) \in \lambda} \frac{1}{(1-q^{h_{\lambda}(i,j) })},\]
where the sum is over all reverse plane partitions $T$ of shape $\lambda$.
\end{cor}

\noindent {\bf Proof}: The coefficient of $q^n$ on the right-hand side is the number of ways to choose $t_{ij} \in \mathbb{N}$ such that $n = \sum_{(i,j)\in\lambda} t_{ij} h_{\lambda}(i,j)$. But by the previous proposition, the $\mathcal{RSK}$ bijection takes such assignments of nonnegative integers to the boxes of $\lambda$ to reverse plane partitions $T$ of shape $\lambda$ whose entries sum to~$n$. The number of such reverse plane partitions is clearly what is counted by the coefficient of $q^n$ on the left-hand side. $\square$

%explain "content"

Actually a careful analysis of our procedure yields a more refined version of this result. Fix some list $\mathbf{x} = (\ldots,x_{-1},x_0,x_1,\ldots)$ of commuting indeterminates indexed by $\mathbb{Z}$; we think of the $x_i$ being nonnegative reals that assign weights to the diagonals of our partition, where $x_0$ corresponds to the main diagonal, $x_1$ the diagonal above the main diagonal, $x_{-1}$ the diagonal below it, and so on. We define the \emph{$\mathbf{x}$-hook-length} of $\lambda$ at a box $(i,j)$ to be
\[ h_{\lambda}(i,j;\mathbf{x}) := \sum_{(i',j') \in H_{\lambda}(i,j)} x_{j'-i'}.\]
Similarly we define the \emph{$\mathbf{x}$-weight} of a reverse plane partition $T = (t_{ij})$ to be
\[ {|T|}_{\mathbf{x}} := \sum_{(i,j) \in \lambda} t_{ij} x_{j-i}.\]
Observe that taking $\mathbf{x} = (\ldots,1,1,1,\ldots)$ recovers the normal hook-lengths and weights. Recall that the \emph{content} of a box $(i,j) \in \lambda$ of a partition is defined to be $j-i$, so these can be thought of as ``content-weighted'' analogues of the classical notions. If $f = \sum_{i}c_ix_i$ where the $c_i \in \mathbb{N}$ and all but finitely many are zero, let us use the notation $\mathbf{z}^f := \prod_{i}z_i^{c_i}$, where the~$z_i$ are some new set of commuting variables also indexed by $\mathbb{Z}$. Exactly the same proofs as above give us the following.

\begin{prop}
Suppose $T \xmapsto{\mathcal{RSK}} \widehat{T}$ with $T = (t_{ij})$  and $\mathrm{sh}(T) = \lambda$. Then
\[{|\widehat{T}|}_{\mathbf{x}} = \sum_{(i,j)\in\lambda} t_{ij} h_{\mathrm{sh}(T)}(i,j;\mathbf{x}). \]
\end{prop}

\begin{cor} \label{cor:weightrppgf}
Let $\lambda$ be a partition. Then
\[ \sum_{T} \mathbf{z}^{|T|_{\mathbf{x}}} = \prod_{(i,j)\in \lambda} \frac{1}{(1-\mathbf{z}^{h_{\lambda}(i,j;\mathbf{x})})},\]
where the sum is over all reverse plane partitions $T$ of shape $\lambda$.
\end{cor}

 The classical generating function for reverse plane partitions in terms of hook-lengths implies the so-called ``hook-length formula'' for standard Young tableaux by taking the limit $q \to 1$ and applying some results from the theory of~$P$-partitions (see~\cite[\S7.22]{stanley2}). Our content-weighted version similarly implies a content-weighted hook-length formula. 
 
 Let $T$ be a standard Young tableau with $\mathrm{sh}(T) \vdash n$. For any $k \in [n]$, let us use $(i_T(k),j_T(k))$ to denote the box whose entry is $k$. Then define
\[ T_{\mathbf{x}} := \prod_{k=1}^{n} \frac{1}{\sum_{l = 1}^{k} x_{j_T(n+1-l) - i_T(n+1-l)}}.\]
\begin{example}
If $T$ is the SYT given by
\[ T = \young(134,25)\]
then we have
\[T_{\mathbf{x}} = \Big(x_0(x_0+x_2)(x_0+x_1+x_2)(x_{-1}+x_0+x_1+x_2)(x_{-1}+2x_0+x_1+x_2)\Big)^{-1}.\]
\end{example}

 \begin{thm} \label{thm:whlf}
 Let $\lambda \vdash n$ be a partition. Then
\[ \sum_{T} T_{\mathbf{x}} = \prod_{(i,j)\in\lambda} \frac{1}{h_{\lambda}(i,j;\mathbf{x})},\]
where the sum is over all standard Young tableaux $T$ of shape $\lambda$.
 \end{thm}
 
\noindent {\bf Proof}: Let $P^{\mathrm{op}}_{\lambda}$ be the order dual to the poset defined in the proof of Proposition~\ref{prop:welldefined}; that is, $P^{\mathrm{op}}_{\lambda}$ is the poset on the boxes of $\lambda$ with~$(i,j) \leq (k,l)$ iff~$i \geq k$ and $j \geq l$. So the upper-left box is maximal in $P^{\mathrm{op}}_{\lambda}$. We recall the main generating function for $P$-partitions~\cite[Theorem 3.15.5]{stanley1}: for any poset~$P = \{p_1,\ldots,p_r\}$ on $r$ elements with a natural labeling $\omega\colon P \to [r]$, we have
 \[ \sum_{\sigma \in \mathcal{A}(P)} y_1^{\sigma(p_1)}\cdots y_r^{\sigma(p_r)} = \sum_{w \in \mathcal{L}(P)} \frac{\prod_{j \in D(w)} y_{w'_1}y_{w'_2}\cdots y_{w'_j}}{\prod_{i=1}^{p} (1-y_{w'_1}y_{w'_2}\cdots y_{w'_i})},\]
 where $\mathcal{A}(P)$ is the set of $P$-partitions, $\mathcal{L}(P)$ is the set of linear extensions of $P$, and $D(w)$ is some set (the descent set of $w$) that shall not concern us. We identify $\mathcal{L}(P)$ with a subset of $\mathfrak{S}_r$, the symmetric group on $r$ letters, by setting
 \[ \mathcal{L}(P) := \{w \in \mathfrak{S}_r\colon w^{-1}\circ\omega\colon P \to [r] \textrm{ is order-preserving}\}.\]
 Then for each $k \in [r]$, we define $w'_k$ to be the integer $j \in [r]$ that satisfies~$w_k = \omega(s_j)$. Applied to our poset $P^{\mathrm{op}}_{\lambda}$, and by substituting $z_{j-i}$ for the variable $y_{u}$ which corresponds to box $u = (i,j)$, this identity amounts to
 \[ \sum_{\substack{T \textrm{ rev. p.p.}, \\ \mathrm{sh}(T) = \lambda}} \mathbf{z}^{|T|_{\mathbf{x}}} = \sum_{\substack{T \textrm{ SYT}, \\ \mathrm{sh}(T) = \lambda}} \frac{m(T)}{\prod_{k=1}^n (1 - \prod_{l=1}^{k}z_{j_T(n+1-l) - i_T(n+1-l)})},\]
where $m(T)$ is some monomial in the $z_i$ which will not concern us. But then by Corollary~\ref{cor:weightrppgf} we obtain
 \[\sum_{\substack{T \textrm{ SYT}, \\ \mathrm{sh}(T) = \lambda}} \frac{m(T)}{\prod_{k=1}^n (1 - \prod_{l=1}^{k}z_{j_T(n+1-l) - i_T(n+1-l)})} = \prod_{(i,j)\in \lambda} \frac{1}{(1-\mathbf{z}^{h_{\lambda}(i,j;\mathbf{x})})}.\]
 Now suppose we set $z_i = 1 - \varepsilon x_i$ for $\varepsilon \in \mathbb{R}_{> 0}$, and we multiply the above equation by $\varepsilon^n$ and then take the limit as $\varepsilon \to 0$. Then for any 
 \[f = (1 - z_{i_1}z_{i_2}\cdots z_{i_k})^{-1},\]
 we have 
 \[\lim_{\varepsilon \to 0} \varepsilon f = (x_{i_1} + x_{i_2} + \cdots + x_{i_k})^{-1};\] 
 in other words, this limit picks out the linear term. Also, for any $m$ which is a monomial in the $z_i$, we clearly get that~$\lim_{\varepsilon \to 0} m = 1$. So in the last equation, the left-hand side goes to
 \[ \sum_{\substack{T \textrm{ SYT}, \\ \mathrm{sh}(T) = \lambda}} \frac{1}{\prod_{k=1}^n \sum_{l=1}^{k} x_{j_T(n+1-l) - i_T(n+1-l)}} =   \sum_{\substack{T \textrm{ SYT}, \\ \mathrm{sh}(T) = \lambda}} T_{\mathbf{x}}\]
 in the limit and the right-hand side goes to $\prod_{(i,j)\in\lambda} \frac{1}{h_{\lambda}(i,j;\mathbf{x})}$, giving us the claimed identity. $\square$
 
 \begin{remark}
Suppose $\lambda \vdash n$. Then note that  setting $\mathbf{x} = (\ldots,1,1,1,\ldots)$ gives $T_{\mathbf{x}} = \frac{1}{n!}$ for any SYT $T$ with $\mathrm{sh}(T) = \lambda$. In this case, Theorem~\ref{thm:whlf} becomes
 \[\sum_{\substack{T \textrm{ SYT},\\ \mathrm{sh}(T) = \lambda}} \frac{1}{n!} = \prod_{(i,j)\in\lambda} \frac{1}{h_{\lambda}(i,j)},\] 
 recovering the usual hook-length formula for the number of SYT of shape~$\lambda$.
 \end{remark}

\section{\texorpdfstring{$\mathcal{RSK}$}{RSK} = RSK: the Greene-Kleitman invariant and the octahedron recurrence} \label{sec:gk}

In this section we will show that the $\mathcal{RSK}$ algorithm defined above via toggles is the same as classical RSK when restricted to matrices. The idea is to show that $\mathcal{RSK}$ has a Greene-Kleitman invariant. To give a rough idea of what this means, Greene's theorem~\cite{greene} says that if we apply RSK to the matrix corresponding to a permutation~$\pi$ and~$\lambda$ is the shape of the resulting SSYTs, then the maximum size of the union of~$k$ disjoint increasing subsequences in $\pi$ is given by~$\lambda_1 + \cdots + \lambda_k$. We want to give a similar description of the output of our $\mathcal{RSK}$ algorithm but where the input matrix, indeed $\mathbb{N}$-tableau, can be arbitrary. In order to establish such a description of the output, it will be helpful to visualize all of the algorithm at once in a three-dimensional array. As it turns out, this array can be described by the famous octahedron recurrence. The connection between RSK and the octahedron recurrence has been noted before by Danilov and Koshevoy (see~\cite{danilov2}, which builds on the work of~\cite{danilov1}).

Let $T = (t_{ij})$ be an $\mathbb{N}$-tableau. We associate to $T$ a three-dimensional array $U_T = (u_{ijk})$ whose coordinates run over all $i,j,k \in \mathbb{Z}$ which satisfy
\begin{enumerate}[(a)]
\item there exists $a,b \in \{0,1,2\}$ such that $(i+a,j+b)$ is a box of $\mathrm{sh}(T)$;
\item $0 \leq k \leq \mathrm{min}(i,j)+1$.
\end{enumerate}
Observe that $U_T$ is ``pyramid-shaped''. For example, if $T$ is a matrix, then~$U_T$ is literally a square pyramid, albeit with the apex above the bottom-right corner. The entries of the array $U_T$ are filled-in as follows: the boundary conditions are defined to be
\[u_{i,j,0} :=  0 \textrm{ and } u_{i,j,\mathrm{min}(i,j)+1} := 0;\]
and for $0 < k < \mathrm{min}(i,j)+1$ we recursively set
\[u_{ijk} := \mathrm{min}(u_{i-1,j,k-1}, u_{i,j-1,k-1}) + \mathrm{max}(u_{i-1,j,k},u_{i,j-1,k}) - u_{i-1,j-1,k-1} + t_{ijk}\]
where
\[ t_{ijk} := \begin{cases} t_{ij} &\textrm{if } k=1 \\
0 &\textrm{otherwise}. \end{cases} \]
The array $U_T$ captures each step of the $\mathcal{RSK}$ algorithm, as Proposition~\ref{prop:uisrsk} makes precise. Two properties of $U_T$ that follow immediately from this definition are that $u_{ijk} = 0$ if $(i,j)$ is not a box of $T$, and that if $S$ is the restriction of $T$ then $U_S$ is the restriction of $U_T$. From $U_T$, we define another, modified three-dimensional array $\overline{U}_T = (\overline{u}_{ijk})$  whose coordinates run over all $i,j,k \in \mathbb{Z}$ which satisfy
\begin{enumerate}[(a)]
\item there exists $a,b \in \{0,1\}$ such that $(i+a,j+b)$ is a box of $\mathrm{sh}(T)$;
\item $0 \leq k \leq \mathrm{min}(i,j)$.
\end{enumerate}
The entries of $\overline{U}_T$ are given by
\[ \overline{u}_{ijk} :=  \sum_{l=0}^{k} u_{ijl}.\]
The intuition behind why we consider the modified array $\overline{U}_T$ is that Greene's theorem concerns partial sums of the form~$\lambda_1 + \cdots + \lambda_k$, so we want to focus on partial sums of our output rather than individual terms. Finally we define a third array $\widetilde{U}_T = (\widetilde{u}_{ijk})$ whose coordinates have the same bounds as $\overline{U}_T$ and whose entries are given by
\[ \widetilde{u}_{ijk} := \overline{u}_{ijk} - \mathrm{rect}_T(i,j).\]
(If $(i,j)$ is not a box of $\mathrm{sh}(T)$, we interpret $\mathrm{rect}_T(i,j)$ as $0$.)  Subtracting the rectangle sum is just a renormalization that allows $\widetilde{U}_T$ to satisfy the octahedron recurrence. We now give an example of these arrays.

\begin{example}
Suppose
\[ A = \begin{pmatrix} 1 & 0 & 2 \\ 0 & 2 & 0 \\ 1 & 1 & 0 \end{pmatrix}\]
as in Example~\ref{ex:rsk}. Then the arrays associated to $A$ are
\begin{center}
$U_A =$ \setlength{\tabcolsep}{9pt}
\begin{tabular}{|c |c | c | c | c | c|}
\hline
$\begin{matrix} 0 & 0 & 0 & 0 &0 \\ 0 & 0 & 0 & 0 &0 \\ 0 & 0 & 0 & 0 &0 \\ 0 & 0 & 0 & 0 &0 \\ 0 & 0 & 0 & 0 &0 \end{matrix}$ & $\begin{matrix} \\ 0 & 0 & 0 & 0 \\ 0 & 1 & 1 & 3 \\ 0 & 1 & 3 &3 \\ 0 & 2 & 4 & 4 \end{matrix}$ & $\begin{matrix} \\ \\ 0 & 0 & 0 \\ 0 & 0 & 2 \\ 0 & 1 & 2 \end{matrix}$ & $\begin{matrix} \\ \\ \\ 0 & 0 \\ 0 & 1 \end{matrix}$ & $\begin{matrix} \\ \\ \\ \\ 0 \end{matrix}$ \\
$k=0$ & $k=1$ & $k=2$ & $k=3$ & $k=4$ \\ \hline
\end{tabular}\\ \medskip

$\overline{U}_A = $\setlength{\tabcolsep}{9pt}
\begin{tabular}{|c | c | c | c | c |}
\hline
$\begin{matrix} 0 & 0 & 0 & 0 \\ 0 & 0 & 0 & 0 \\ 0 & 0 & 0 & 0 \\ 0 & 0 & 0 & 0 \end{matrix}$ & $\begin{matrix} \\ 1 & 1 & 3 \\  1 & 3 & 3 \\  2 & 4 & 4 \end{matrix}$ & $\begin{matrix}  \\ \\ 3 & 5 \\ 5 & 6 \end{matrix}$ & $\begin{matrix} \\ \\ \\ 7 \end{matrix}$ \\
$k=0$ & $k=1$ & $k=2$ & $k=3$ \\ \hline
\end{tabular} \\ \medskip

$\widetilde{U}_A = $\setlength{\tabcolsep}{9pt}
\begin{tabular}{|c | c | c | c | c |}
\hline
$\begin{matrix} 0 & 0 & 0 & 0 \\ 0 & -1 & -1 & -3 \\ 0 & -1 & -3 & -5 \\ 0 & -2 & -5 & -7 \end{matrix}$ & $\begin{matrix} \\ 0 & 0 & 0 \\  0 & 0 & -2 \\  0 & -1 & -3 \end{matrix}$ & $\begin{matrix}  \\ \\ 0 & 0 \\ 0 & -1 \end{matrix}$ & $\begin{matrix} \\ \\ \\ 0 \end{matrix}$ \\
$k=0$ & $k=1$ & $k=2$ & $k=3$ \\ \hline
\end{tabular}
\end{center}
Here the bottom-right corner of each level in each array is aligned and this corner has $(i,j) = (3,3)$. Within each level we use normal matrix coordinates so that the entry with minimal $(i,j)$ is the upper-left corner. $\square$
\end{example}

\begin{prop} \label{prop:uisrsk}
Suppose $T \xmapsto{\mathcal{RSK}} \widehat{T}$ where $ \widehat{T} = (\widehat{t}_{ij})$. Let $(i,j) \in \mathrm{sh}(T)$, and suppose $m$ is the unique natural number such that $(i+m,j+m)$ is a border box of $\lambda$. Let $U_T = (u_{ijk})$ be as defined above. Then~$\widehat{t}_{ij} = u_{i+m,j+m,m+1}$.
\end{prop}

\noindent {\bf Proof}: This is clear from construction. It can be easily verified by again considering inductively the addition of a box. $\square$

\begin{prop}
Let $T = (t_{ij})$ be an $\mathbb{N}$-tableau. Then the array $\widetilde{U}_T = (\widetilde{u}_{ijk})$ defined above can also be described as follows: the boundary conditions are given by $\widetilde{u}_{i,j,0} = -\mathrm{rect}_T(i,j)$ and $\widetilde{u}_{i,j,\mathrm{min}(i,j)} = 0$; and for $0 < k < \mathrm{min}(i,j)$ we have the recursive formula
\[ \widetilde{u}_{ijk} = \mathrm{max}(\widetilde{u}_{i-1,j,k} + \widetilde{u}_{i,j-1,k-1},\widetilde{u}_{i-1,j,k-1} + \widetilde{u}_{i,j-1,k}) - \widetilde{u}_{i-1,j-1,k-1}.\]
In other words, $\widetilde{U}_T$ satisfies the tropical octahedron recurrence.
\end{prop}

\noindent {\bf Proof}: The boundary condition~$\widetilde{u}_{i,j,0} = -\mathrm{rect}_T(i,j)$ is satisfied by definition. To verify the boundary condition $\widetilde{u}_{i,j,\mathrm{min}(i,j)} = 0$ requires a little more work. Suppose $U_T$ is as defined above and also that $T \xmapsto{\mathcal{RSK}} \widehat{T}$.  By induction we may assume that this boundary condition holds for smaller tableaux. In particular, suppose $T$ is obtained from another tableau by adding a corner box $(i,j)$ and the boundary condition holds for this smaller tableau. Then clearly $\widetilde{u}_{i',j',\mathrm{min}(i',j')} = 0$ for all $(i',j') \neq (i,j)$ because this matrix entry will be the same as in the one corresponding to the smaller tableau. But also, Proposition~\ref{prop:uisrsk} tells us that~$\sum_{l=0}^{\mathrm{min}(i,j)} u_{ijl} = \mathrm{diag}_{\widehat{T}}(i,j)$. And Proposition~\ref{prop:diagsum} tells us that $\mathrm{diag}_{\widehat{T}}(i,j) = \mathrm{rect}_T(i,j)$. So indeed we can conclude~$\widetilde{u}_{i,j,\mathrm{min}(i,j)} = 0$ as well, and the condition follows by induction.

Now we need to check the recursive formula. Let $\overline{U}_T$ be as defined above and let $0 < k < \mathrm{min}(i,j)$; we first claim that
\[ \overline{u}_{ijk} = \mathrm{max}(\overline{u}_{i-1,j,k} + \overline{u}_{i,j-1,k-1},  \overline{u}_{i-1,j,k-1} + \overline{u}_{i,j-1,k}) - \overline{u}_{i-1,j-1,k-1} + t_{ij}. \tag{*} \label{eqn:overu} \]
We prove~(\ref{eqn:overu}) by induction on $k$. The case $k=1$ holds because 
\begin{align*}
u_{ij1} &= \mathrm{min}(u_{i-1,j,0},u_{i,j-1,0}) + \mathrm{max}(u_{i-1,j,k},u_{i,j-1,k}) - u_{i-1,j-1,0} + t_{i,j} \\
&= \mathrm{max}(u_{i-1,j,k},u_{i,j-1,k}) + t_{i,j}.
\end{align*}
So now assume $k > 1$ and the result is known for smaller $k$. Then the induction hypothesis and a routine computation give
\begin{align*}
\overline{u}_{ijk} = &\overline{u}_{i,j,k-1} + u_{ijk} \\
= &\mathrm{max}(\overline{u}_{i-1,j,k-1} + \overline{u}_{i,j-1,k-2},  \overline{u}_{i-1,j,k-2} + \overline{u}_{i,j-1,k-1}) - \overline{u}_{i-1,j-1,k-2} + t_{ij} \\
&+\mathrm{min}(u_{i-1,j,k-1},u_{i,j-1,k-1}) + \mathrm{max}(u_{i-1,j,k}, u_{i,j-1,k}) - u_{i-1,j-1,k-1} \\
= & \overline{u}_{i-1,j,k-2} + \overline{u}_{i,j-1,k-2} + \mathrm{max}(u_{i-1,j,k-1},u_{i,j-1,k-1})  - \overline{u}_{i-1,j-1,k-1} \\
 &+ \mathrm{min}(u_{i-1,j,k-1},u_{i,j-1,k-1}) + \mathrm{max}(u_{i-1,j,k}, u_{i,j-1,k})  + t_{ij} \\
 = & \overline{u}_{i-1,j,k-2} + \overline{u}_{i,j-1,k-2} + u_{i-1,j,k-1} + u_{i,j-1,k-1} \\
 &+ \mathrm{max}(u_{i-1,j,k}, u_{i,j-1,k})  - \overline{u}_{i-1,j-1,k-1} + t_{ij}\\
 =&\overline{u}_{i-1,j,k-1} + \overline{u}_{i,j-1,k-1} + \mathrm{max}(u_{i-1,j,k}, u_{i,j-1,k})  - \overline{u}_{i-1,j-1,k-1} + t_{ij} \\
 =&\mathrm{max}(\overline{u}_{i-1,j,k} + \overline{u}_{i,j-1,k-1},  \overline{u}_{i-1,j,k-1} + \overline{u}_{i,j-1,k}) - \overline{u}_{i-1,j-1,k-1} + t_{ij}.
\end{align*}
Thus by induction~(\ref{eqn:overu}) holds for all $k$. Finally for $0 < k < \mathrm{min}(i,j)$ we compute
\begin{align*}
\widetilde{u}_{ijk} = &-\mathrm{rect}_T(i,j) + \overline{u}_{ijk} \\
= &-\mathrm{rect}_T(i,j) + t_{ij} \\
&+ \mathrm{max}(\overline{u}_{i-1,j,k} + \overline{u}_{i,j-1,k-1},  \overline{u}_{i-1,j,k-1} + \overline{u}_{i,j-1,k}) - \overline{u}_{i-1,j-1,k-1}  \\
= &-\mathrm{rect}_{T}(i-1,j) - \mathrm{rect}_T(i,j-1) + \mathrm{rect}_T(i-1,j-1) \\
&+ \mathrm{max}(\overline{u}_{i-1,j,k} + \overline{u}_{i,j-1,k-1},  \overline{u}_{i-1,j,k-1} + \overline{u}_{i,j-1,k}) - \overline{u}_{i-1,j-1,k-1}  \\
= &  \mathrm{max}(\widetilde{u}_{i-1,j,k} + \widetilde{u}_{i,j-1,k-1},\widetilde{u}_{i-1,j,k-1} + \widetilde{u}_{i,j-1,k}) - \widetilde{u}_{i-1,j-1,k-1},
\end{align*}
and the proposition is proved.  $\square$

We now state the Greene-Kleitman invariant for $\mathcal{RSK}$, which requires a bit of terminology.

\begin{definition}
Let $\lambda$ be a partition and $(a,b), (c,d)$ boxes of $\lambda$, with $a \leq c$ and $b \leq d$. A \emph{path} $p$ from $(a,b)$ to $(c,d)$ in $\lambda$ is a sequence
\[ (a,b) = (i_0,j_0), (i_1,j_1), \ldots, (i_l,j_l) = (c,d)\]
of boxes in $\lambda$, where for $1 \leq k \leq l$ we have either $i_k = i_{k-1} + 1$ and $j_k = j_{k-1}$, or $i_k = i_{k-1}$ and $j_k = j_{k-1} + 1$. We use the notation
\[ \mathrm{box}(p) := \{ (i_k,j_k)\colon k = 0,\ldots,l\}\]
to denote the set of boxes in $p$. Let $T = (t_{ij})$ be an $\mathbb{N}$-tableau of shape $\lambda$. The weight of $p$ in $T$ is then defined to be $\mathrm{wt}_T(p) := \sum_{(i,j) \in \mathrm{box}(p)} t_{i,j}$. Two paths $p$ and $q$ in $\lambda$ are \emph{noncrossing} if $\mathrm{box}(p) \cap \mathrm{box}(q) = \emptyset$. We say that several paths $p_1,\ldots,p_k$ are noncrossing if they are pairwise noncrossing. Suppose that~$\mathcal{X} = \{x_i\}_{i=1}^{k}$ and $\mathcal{Y} = \{y_i\}_{i=1}^{k}$ are sequences of boxes in $\lambda$ (with the same number of terms). Then we define
\[ \mathrm{NCPath}(\mathcal{X},\mathcal{Y}) := \left\{ \{p_1,\ldots,p_k\} \colon \parbox{2.5in}{\begin{center}$p_1,\ldots,p_k$ are noncrossing paths in $\lambda$ \\ and $p_i$ is a path from $x_i$ to $y_i$ \\ for all $1 \leq i \leq k$ \end{center}}\right\} \]
to be the set of noncrossing paths connecting $\mathcal{X}$ and $\mathcal{Y}$.
\end{definition}

\begin{thm} \label{thm:gk}
Let $T = (t_{ij})$ be an $\mathbb{N}$-tableau of shape $\lambda$. Let $(i,j)$ be a box of~$\lambda$ and $1 \leq k \leq \mathrm{min}(i,j)$. Set $\mathcal{X} := \{ (1,l)\}_{l=1}^{k}$ and $\mathcal{Y} := \{ (i,j-k+l)\}_{l=1}^{k}$ and define
\[ m(i,j,k) := \mathrm{max}\left \{ \sum_{l=1}^{k} \mathrm{wt}_T(p_l)\colon \{p_1,\ldots,p_k\} \in \mathrm{NCPath}(\mathcal{X},\mathcal{Y}) \right\}.\]
Let $\overline{U}_T = (\overline{u}_{ijk})$ be as defined above. Then $m(i,j,k) = \overline{u}_{ijk}$.
\end{thm}

The key to proving Theorem~\ref{thm:gk}, which we will not do here, is to show that (the appropriate normalization of) the array of values $m(i,j,k)$ satisfies the tropical octahedron recurrence. The birational analogue of this fact is~\cite[Theorem 1]{danilov3}. Of course any birational identity tropicalizes and so that result implies the combinatorial version of the theorem (i.e., the version as stated). A bijective proof of the birational analogue is also given in~\cite{farber}. This Greene-Kleitman invariant for $\mathcal{RSK}$ implies that, restricted to matrices, $\mathcal{RSK}$ is the same as classical RSK because RSK respects the same invariant. However, the fact that RSK has such a Greene-Kleitman invariant is apparently folklore; the best reference we have for it is~\cite[Theorem~12]{krattenthaler}.

\section{Future directions: dual RSK; other \texorpdfstring{$d$}{d}-complete posets} \label{sec:future}

We briefly describe two different threads of possible future research:

\begin{enumerate}[(1)]
\item There is another correspondence involving matrices and SSYT, described by Knuth in the same seminal paper where he introduced RSK~\cite{knuth}, called dual RSK. Whereas RSK is a bijection between $\mathbb{N}$-matrices and pairs of SSYT of the same shape, dual RSK is a bijection between $(0,1)$-matrices and pairs of SSYT whose shapes are conjugate. Normal RSK yields the Cauchy identity
\[ \sum_{\lambda} s_{\lambda}(x_1,x_2,\ldots) s_{\lambda}(y_1,y_2,\ldots) = \prod_{i,j \geq 1} \frac{1}{(1-x_iy_j)}\]
for Schur functions $s_{\lambda}$, while dual RSK yields the dual Cauchy identity
\[ \sum_{\lambda} s_{\lambda}(x_1,x_2,\ldots) s_{\lambda'}(y_1,y_2,\ldots) = \prod_{i,j \geq 1} (1+x_iy_j).\]
See for example~\cite[\S 7.14]{stanley2}. Another way to think about the duality is that normal RSK is bosonic while dual RSK is fermionic. It would be interesting to extend the above approach, using local toggle operations, to dual RSK as well. Although birational dual RSK has already been studied by Noumi and Yamada~\cite{noumi}, their account is strictly algebraic and does not give a combinatorial description of what is going on in the dual case. Most interesting would be if birational RSK and birational dual RSK can be realized simultaneously in some kind of superalgebra with both commuting (bosonic) and anticommuting (fermionic) variables.

\item Proctor~\cite{proctor} has defined $d$-complete posets which generalize Young diagrams. In particular, these posets have a hook-length formula enumerating their number of linear extensions. It would be very interesting to try to extend the toggle operation described here on Young diagrams to $d$-complete posets and in the process obtain a bijective proof of the hook-length formula for these more general objects. To our knowledge, no such bijective proof exists.

\end{enumerate}

\bibliography{rsk}{}
\bibliographystyle{plain}

\end{document}